\documentclass[11pt]{article}
\usepackage{amsfonts}
\usepackage{bbm}
\usepackage{mathrsfs}
\usepackage{amscd}
\usepackage{color}
\usepackage{amsmath,amsfonts,amssymb,amscd}
\usepackage{indentfirst,graphics,epsfig,psfrag}
\input{epsf}
\usepackage{ifpdf}
\usepackage{enumerate}
\usepackage{appendix}
\usepackage{enumerate}

\usepackage{lineno}

\makeatletter \@addtoreset{figure}{section} \makeatother
\makeatletter
\long\def\@makecaption#1#2{%
   \vskip 10\p@
   \setbox\@tempboxa\hbox{{#1}\ \ #2}%
   \ifdim \wd\@tempboxa >\hsize

       {#1}\ \ #2\par
   \else
       \hbox to\hsize{\hfil\box\@tempboxa\hfil}%
   \fi}
\makeatother

\begin{document}

\title{Resistance distance and Kirchhoff index in the corona-vertex and the
corona $-$ edge of subdivision graph
\thanks{\footnotesize }}
\author{Qun Liu~$^{a}$, Jia-Bao Liu~$^{b,\dag}$, Shaohui Wang~$^c$, Masoud Karimi~$^{d,e}$\setcounter{footnote}{-1 }\thanks{${~\dag}$ Corresponding
author. E-mail addresses:
liuqun09@yeah.net, liujiabaoad@163.com, shaohuiwang@yahoo.com, karimimth@yahoo.com.}\\
      \emph{\small a. School of Mathematics and Statistics, Hexi University,}\\
   \emph{\small Gansu, Zhangye, 734000, P.R. China}\\
      \emph{\small b. School of Mathematic and Physics, Anhui Jianzhu University,}\\
   \emph{\small Hefei, 230601, P.R. China}\\
 \emph{\small c. Department of Mathematic and Computer Science, Adelphi University,}\\
   \emph{\small Garden City, NY 11550, USA}\\
    \emph{\small d. School of Mathematical Sciences, Anhui University,}\\
   \emph{\small Hefei, 230601, P.R. China}\\
   \emph{\small e. Department of Mathematics, Bojnourd Branch, Islamic Azad University,}\\
   \emph{\small Bojnourd, Iran}\\}
\title{\bf{Resistance distance and Kirchhoff index in corona-vertex and
corona-edge of subdivision graph}} \frenchspacing
\date{}
\maketitle

\begin{minipage}{13cm}

{\bf Abstract}


\vskip 0.05in

{\small   The subdivision graph $S(G)$ of a graph $G$ is the graph obtained by inserting a new vertex into every edge of $G$. In $\cite{PL}$, two classes of new corona graphs,
the corona-vertex of the subdivision graph
$G_{1}\diamondsuit G_{2}$ and corona-edge of the subdivision
graph $G_{1}\star G_{2}$ were defined. The adjacency spectrum
and the signless Laplacian spectrum of the two new graphs
were computed when $G_{1}$ is an arbitrary graph and $G_{2}$
is an $r$-regular graph. In this paper, we give the formulate of the resistance distance and
the Kirchhoff index in $G_{1}\diamondsuit G_{2}$
and $G_{1}\star G_{2}$ when $G_{1}$ and $G_{2}$
are arbitrary graphs. These results generalize
them in $\cite{PL}$.
}

\vskip 0.1in

\end{minipage}

\vskip 0.1in


\section{ Introduction }
Throughout this article, all graphs considered are simple and
undirected. Let $G=(V(G),$
$E(G))$ be a graph with vertex set
$V(G)=\{v_{1},v_{2},...,v_{n}\}$ and edge set
$E(G)=\{e_{1},e_{2},...,e_{m}\}$.
In 1993, Klein and Randi$\acute{c}$ $\cite{DJK}$ introduced a distance function named resistance
distance on the basis of electrical network theory. They view a graph
as an electrical network and each edge of the graph is assumed to be a unit resistor,
then take the  resistance distance between vertices to be the effective resistance
between them. Let $G$ be a simple graph with the vertex set $V(G)=\{v_{1},v_{2},...,v_{n}\}$,
and $r_{ij}(G)$ denote the effective resistance distance between vertices $v_{i}$ and $v_{j}$
as computed with Ohm${'}$s law when all the edges of $G$ are considered to be unit resistors.
The sum of resistance distance $Kf(G)=\sum_{u,v\in V(G)}r_{uv}(G)$ was proposed in Ref. $\cite{DJK}$, later called the Kirchhoff index of $G$ in
 Ref. $\cite{DB}$. The resistance distance
and the Kirchhoff index attracted extensive attention due to its wide applications in physics, chemistry and others.
For more information on resistance distance and Kirchhoff index of graphs, the readers are referred to Refs. 
[8-10, 14-19, 21-29] and the references therein.

Let $G=(V(G),E(G))$ be a graph with vertex set $V(G)$ and edge set $E(G)$.
Let $d_{i}$ be the degree of vertex $i$ in $G$ and $D_{G}=diag(d_{1},
d_{2},\cdots, d_{|V(G)|})$ the diagonal matrix with all vertex degrees of
$G$ as its diagonal entries. For a graph $G$, let $A_{G}$ and $R_{G}$ denote the adjacency matrix
and vertex-edge incidence matrix of $G$, respectively. The matrix
$L_{G}=D_{G}-A_{G}$ is called the Laplacian matrix of $G$, where
$D_{G}$ is the diagonal matrix of vertex degrees of $G$. The line graph $l(G)$ of
a graph $G$ is the graph $V(l(G))=E(G)$ such that two vertices in $l(G)$
are adjacent if and only if the corresponding edges in $G$ are adjacent.
Note that
$R_{G}R_{G}^{T}=D_{G}+A_{G}$ and $R_{G}^{T}R_{G}=A(l(G))+2I_{m}$, where
$m=|E(G)|$.
The $\{1\}$-inverse of $M$ is a matrix $X$ such that $MXM=M$. If $M$ is singular, then it has infinite
$\{1\}$- inverse $\cite{AB}$. We use $M^{(1)}$ to denote any $\{{1}\}$- inverse of a matrix $M$,
and let $(M)_{uv}$ denote the $(u,v)$- entry of $M$. It is known that resistance distances in a connected graph $G$ can be obtained from any $\{1\}$- inverse of $L_G$ $(\cite{RB},\cite{RBB})$.

The subdivision graph $S(G)$ of a graph $G$ is the graph obtained
by inserting a new vertex into every edge of $G$. The set of such
new vertices is denoted by $I(G)$. In Ref. $\cite{PL}$, two new graph
operations based on subdivision graphs: the corona-vertex
and the corona-edge of the subdivision graph were introduced.
The corona-vertex of the subdivision
graph of two vertex-disjoint
graphs $G_{1}$ and $G_{2}$ denoted by $G_{1}\diamondsuit G_{2}$,
is the graph obtained from $G_{1}$ and $|V(G_{1})|$ copies of
$S(G_{2})$ by joining the ith vertex of $V(G_{1})$ to every
vertex in the ith copy of $G_{2}$.
The corona-edge of the subdivision graph of two vertex-disjoint
graphs $G_{1}$ and $G_{2}$, denoted by $G_{1}\star G_{2}$,
is the graph obtained from $G_{1}$ and $|V(G_{1})|$ copies of
$S(G_{2})$ by joining the ith vertex of $V(G_{1})$ to every
vertex in the ith copy of $I(G_{2})$.

Note that if $G_{1}$ is a graph on $n_{1}$ vertices and $m_{1}$ edges,
and $G_{2}$ is a graph on $n_{2}$ vertices and $m_{2}$ edges, then the
corona-vertex of the subdivision graph $G_{1}\diamondsuit G_{2}$ has
$n_{1}(1+n_{2}+m_{2})$ vertices and $m_{1}+n_{1}n_{2}+2n_{1}m_{2}$
edges, and the corona-edge of the subdivision graph $G_{1}\star G_{2}$
has $n_{1}(1+n_{2}+m_{2})$ vertices and $m_{1}+3n_{1}m_{2}$ edges.

It is of interest to study the Kirchhoff index of graphs derived from a single graph. The formulas
and lower bounds of the Kirchhoff index of the line graph, subdivision graph
of a connected regular graph are reported in Refs. [5,18]. The formulas for the
resistance distance and Kirchhoff index of the subdivision graphs of general graphs were obtained in
Refs. $(\cite{YJ},\cite{YJY})$. In Ref. $\cite{CB}$, Bu et al. obtained the formulae for resistance distances in subdivision-vertex
join and subdivision-edge join. In $\cite{XL}$, Liu gave the resistance distance and Kirchhoff index
of $R$-vertex join and $R$-edge join of two graphs.
Recently Liu et al $\cite{LiuPF2016}$ gave the $\{1\}$-inverse of the Laplacian of subdivision-vertex
and subdivision-edge coronae.
Motivated by these work, in this paper, we formulate the resistance distance and the Kirchhoff index in $G_{1}\diamondsuit G_{2}$
and $G_{1}\star G_{2}$ respectively.

\section{Preliminaries}
For a square matrix $M$, the group inverse of $M$, denoted by $M^{\#}$, is the unique matrix $X$ such that
$MXM=M$, $XMX=X$ and $MX=XM$. It is known that $M^{\#}$ exists if and only if $rank(M)=rank(M^{2})$. If $M$
is real symmetric, then $M^{\#}$ exists and $M^{\#}$ is a symmetric $\{1\}$- inverse of $M$. Actually,
$M^{\#}$ is equal to the Moore-Penrose inverse of $M$ since $M$ is symmetric.

We use $M^{(1)}$ to denote any $\{1\}$- inverse of a matrix $M$. Let $(M)_{uv}$ denote the $(u,v)$ entry of $M$.

{\bf Lemma 2.1}$(\cite{RB},\cite{CLJ}$) \ Let $G$ be a connected graph. Then
$$r_{uv}(G)=(L^{(1)}_{G})_{uu}+(L^{(1)}_{G})_{vv}-(L^{(1)}_{G})_{uv}-(L^{(1)}_{G})_{vu}=(L^{\#}_{G})_{uu}+(L^{\#}_{G})_{vv}-2(L^{\#}_{G})_{uv}.$$

Let $\mathbf1_{n}$ denotes the column vector of dimension $n$ with all the entries equal one.
We will often use $\mathbf{1}$ to denote an all-ones column vector if the dimension can be read from the
context.
\vskip 0.1in
{\bf Lemma 2.2} $(\cite{LZ})$ \ For any graph, we have
$L^{\#}_{G}\mathbf{1}=0$.

\vskip 0.1in
{\bf Lemma 2.3} $(\cite{JZ})$ \ Let
 \[
\begin{array}{crl}
M=\left(
  \begin{array}{cccccccccccccccc}
   A& B  \\
   C & D \\
 \end{array}
  \right)
\end{array}
\]
be a nonsingular matrix. If $A$ and $D$ are nonsingular,
then
\[
\begin{array}{crl}
M^{-1}&=&\left(
  \begin{array}{cccccccccccccccc}
   A^{-1}+A^{-1}BS^{-1}CA^{-1}&-A^{-1}BS^{-1} \\
   -S^{-1}CA^{-1} & S^{-1}\\
 \end{array}
  \right)
\\&=&\left(
  \begin{array}{cccccccccccccccc}
   (A-BD^{-1}C)^{-1} &-A^{-1}BS^{-1} \\
   -S^{-1}CA^{-1} & S^{-1}\\
 \end{array}
  \right),
\end{array}
\]
where $S=D-CA^{-1}B.$

{\bf Lemma 2.4} $(\cite{LZ})$ \ Let
 \[
\begin{array}{crl}
M=\left(
  \begin{array}{cccccccccccccccc}
   L_{1}& L_{2} \\
   L_{2}^{T} & L_{3} \\
 \end{array}
  \right)
\end{array}
\]
be a symmetric  matrix and $L_{1}$ is nonsingular.
Then
\[
\begin{array}{crl}
X&=&\left(
  \begin{array}{cccccccccccccccc}
   L_{1}^{-1}+L_{1}^{-1}L_{2}S^{\#}L_{2} ^{T}L_{1}^{-1}&-L_{1}^{-1}L_{2}S^{\#} \\
   -S^{\#}L_{2}^{T}L_{1}{-1} & S^{\#}\\
 \end{array}
  \right)
\end{array}\frac{}{}
\]is a symmetric $\{1\}$-inverse of $M$,
where $S=L_{3} -L_{2}^{T}L_{1}^{-1}L_{2}.$
\vskip 0.1in
For a vertex $i$ of a graph $G$, let $T(i)$ denote the set of all neighbors of $i$ in $G$.
\vskip 0.1in
{\bf Lemma 2.5}$(\cite{CB})$ \ Let $G$ be a connected graph. For any $i,j\in V(G)$,
$$r_{ij}(G)=d^{-1}_{i}(1+\sum_{k\in T(i)}r_{kj}(G)-d^{-1}_{i}\sum_{k,l\in T(i)}r_{kl}(G)).$$

\vskip 0.1in
{\bf Lemma 2.6}$(\cite{LZ})$ \ Let $G$ be a connected graph on $n$ vertices. Then
$$Kf(G)=ntr(L^{(1)}_{G})-\mathbf{1}^{T}L^{(1)}_{G}\mathbf{1}=ntr(L^{\#}_{G}).$$

{\bf Lemma 2.7} \ Let $\mathbf1_{n}$ and $j_{n\times m }$
be all-one column vector of dimensions $n$ and all-one
$n\times m$ matrix, respectively. For any matrix $A$, we have
$$(\mathbf1_{m_{2}}\otimes I_{n_{1}})A_{n_{1}\times n_{1}}(\mathbf1^{T}_{n_{2}}\otimes I_{n_{1}})=
j_{m_{2}\times n_{2}}\otimes A_{n_{1}\times n_{1}},$$
where $\otimes$ denote the Kronecker product of matrices $A=(a_{ij})$ and $B$, i.e $A\otimes B$ is
defined to be the partition matrix $(a_{ij}B)$.

{\bf Proof.} \
\[
\begin{array}{crl}
(\mathbf1_{m_{2}}\otimes I_{n_{1}})A_{n_{1}\times n_{1}}(\mathbf1^{T}_{n_{2}}\otimes I_{n_{1}})
&=&\left(
  \begin{array}{cccccccccccccccc}
    I_{n_{1}}\\
  I_{n_{1}} \\
 \vdots\\
  I_{n_{1}} \\
 \end{array}
  \right)_{(n_{1}m_{2})\times n_{1}}A_{n_{1}\times n_{1}}
 \begin{array}{crl}
 \left(
  \begin{array}{cccccccccccccccc}
    I_{n_{1}}&I_{n_{1}}&\cdots &I_{n_{1}}\\
 \end{array}
 \right)_{n_{1}\times (n_{1}n_{2})}
\end{array}\\
&=&\left(
  \begin{array}{cccccccccccccccc}
  A&A&\cdots& A\\
  A&A&\cdots &A\\
  \vdots&\vdots&\ddots&\vdots \\
  A&A&\cdots &A\\
 \end{array}
 \right)_{(n_{1}m_{2})\times (n_{1}n_{2})}\\
 &=&j_{m_{2}\times n_{2}}\otimes A_{n_{1}\times n_{1}}.
\end{array}
\]
\hfill$\square$

\section{Resistance distance in corona-vertex and
corona-edge of subdivision graph}
We first give formulae for resistance distance in $G_{1}\diamondsuit G_{2}$.
\vskip 0.1in

{\bf Theorem 3.1} \ Let $G_{1}$ be a graph on $n_{1}$ vertices and $m_{1}$
edges and $G_{2}$ be a graph on $n_{2}$ vertices and $m_{2}$ edges.
Then the following hold:

(i) For any $i,j\in V(G_{1})$, we have
\begin{eqnarray*}
r_{ij}(G_{1}\diamondsuit G_{2})&=&(L_{G_{1}}^{\#})_{ii}+(L_{G_{1}}^{\#})_{jj}-2(L_{G_{1}}^{\#})_{ij}.
\end{eqnarray*}

(ii) For any $i,j\in V(G_{2})$, we have

~~~~~~~~~$r_{ij}(G_{1} \diamondsuit G_{2})=2((L_{G_{2}}+2I_{n_{2}})^{-1}\otimes I_{n_{1}})_{ii}+2((L_{G_{2}}+2I_{n_{2}})^{-1}\otimes I_{n_{1}})_{jj}$
\vskip 0.1in
$-2((L_{G_{2}}+2I_{n_{2}})^{-1}\otimes I_{n_{1}})_{ij}.$

(iii) For any $i\in V(G_{1})$, $j\in V(G_{2})$, we have
\begin{eqnarray*}
r_{ij}(G_{1}\diamondsuit G_{2})&=&(L_{G_{1}}^{\#})_{ii}+2((L_{G_{2}}+2I_{n_{2}})^{-1}\otimes I_{n_{1}})_{jj}-(L_{G_{1}}^{\#})_{ij}.
\end{eqnarray*}

(iv) For any $i\in I(G_{2})$, $j\in V(G_{1})\cup V(G_{2})$, let $u_{i}v_{i}\in E(G_{2})$ denote the edge corresponding
to $i$, we have
\begin{eqnarray*}
r_{ij}(G_{1}\diamondsuit G_{2})&=&\frac{1}{2}+\frac{1}{2}r_{u_{i}j}(G_{1}\diamondsuit G_{2})+\frac{1}{2}r_{v_{i}j}(G_{1}\diamondsuit G_{2})-\frac{1}{4}r_{u_{i}v_{i}}(G_{1}\diamondsuit G_{2}).
\end{eqnarray*}

(v) For any $i,j\in I(G_{2})$, let $u_{i}v_{i}, u_{j}v_{j}\in E(G_{2})$ denote the edges corresponding
to $i,j$, we have
\begin{eqnarray*}
r_{ij}(G_{1}\diamondsuit G_{2})
&=&\frac{1}{2}+\frac{1}{2}r_{u_{i}j}(G_{1}\diamondsuit G_{2})+\frac{1}{2}r_{v_{i}j}(G_{1}\diamondsuit G_{2})-\frac{1}{4}r_{u_{i}v_{i}}(G_{1}\diamondsuit G_{2})\\
&=& 1+\frac{1}{4}[r_{u_{i}u_{j}}(G_{1}\diamondsuit G_{2})+r_{u_{i}v_{j}}(G_{1}\diamondsuit G_{2})+r_{v_{i}u_{j}}(G_{1}\diamondsuit G_{2})\\
&&+r_{v_{i}v_{j}}(G_{1}\diamondsuit G_{2})-r_{u_{i}v_{i}}(G_{1}\diamondsuit G_{2})-r_{u_{j}v_{j}}(G_{1}\diamondsuit G_{2})].
\end{eqnarray*}

{\bf Proof.} \ Let $R_{2}$ be the incidence matrix of $G_{2}$. Then with
a proper labeling of vertices (see Ref.$\cite{PL}$),
the Laplacian matrix of $G_{1}\diamondsuit G_{2}$ can be written as
\[
\begin{array}{crl}
 L(G_{1}\diamondsuit G_{2})=\left(
  \begin{array}{cccccccccccccccc}
    2I_{m_{2}}\otimes I_{n_{1}}& -R_{2}^{T}\otimes I_{n_{1}}&0_{n_{1}m_{2}\times n_{1}}\\
   -R_{2}\otimes I_{n_{1}} &(D_{G_{2}}+I_{n_{2}})\otimes I_{n_{1}}&-\mathbf1_{n_{2}}\otimes I_{n_{1}}\\
   0_{n_{1}\times n_{1}m_{2}}&-\mathbf1^{T}_{n_{2}}\otimes I_{n_{1}}&L_{G_{1}}+n_{2}I_{n_{1}}\\
 \end{array}
  \right)
\end{array},
\]
where $0_{s,t}$ denotes the $s\times t$ matrix with all entries equal to zero. Let
\[
\begin{array}{crl}
 L_{1}=\left(
  \begin{array}{cccccccccccccccc}
    2I_{m_{2}}\otimes I_{n_{1}}& -R_{2}^{T}\otimes I_{n_{1}}\\
   -R_{2}\otimes I_{n_{1}} &(D_{G_{2}}+I_{n_{2}})\otimes I_{n_{1}}\\
 \end{array}
  \right)
\end{array},
\]
 $L_{2}=\left(
  \begin{array}{cc}
    0_{n_{1}m_{2}\times n_{1}}\\
   -\mathbf1_{n_{2}}\otimes I_{n_{1}} \\
  \end{array}
\right)$, $L_{2}^{T}=\left(
  \begin{array}{cc}
   0_{n_{1}\times n_{1}m_{2}}&-\mathbf1^{T}_{n_{2}}\otimes I_{n_{1}}\\
  \end{array}
\right)$, and $L_{3}=L_{G_{1}}+n_{2}I_{n_{1}}$.

We begin with the calculation about $L_{1}^{-1}$. By Lemma 2.3, we have
\begin{eqnarray*}
R_{1}&=&[(D_{G_{2}}+I_{n_{2}})\otimes I_{n_{1}}-(R_{2}\otimes I_{n_{1}})(2I_{m_{2}}\otimes I_{n_{1}})^{-1}(R_{2}^{T}\otimes I_{n_{1}})]^{-1}\\
&=&[(D_{G_{2}}+I_{n_{2}})\otimes I_{n_{1}}-\frac{1}{2}R_{2}R^{T}_{2}\otimes I_{n_{1}}]^{-1}\\
&=&[(D_{G_{2}}+I_{n_{2}})-\frac{1}{2}(D_{G_{2}}+A_{G_{2}})]^{-1}\otimes I_{n_{1}}\\
&=&2(L_{G_{2}}+2I_{n_{2}})^{-1}\otimes I_{n_{1}}.
\end{eqnarray*}
One can obtain that
\[
\begin{array}{crl}
L_{1}^{-1}=\left(
  \begin{array}{cccccccccccccccc}
    T_{1}& (R^{T}_{2}(L_{G_{2}}+2I_{n_{2}})^{-1})\otimes I_{n_{1}}\\
  ((L_{G_{2}}+2I_{n_{2}})^{-1}R_{2})\otimes I_{n_{1}} & 2(L_{G_{2}}+2I_{n_{2}})^{-1}\otimes I_{n_{1}}\\
 \end{array}
  \right),
\end{array}
\]
where $T_{1}=\frac{1}{2}(I_{m_{2}}+R^{T}_{2}(2I_{n_{2}}+L_{G_{2}})^{-1}R_{2})\otimes I_{n_{1}}$.
\vskip 0.1in
Now we are ready to calculate $S$.
Let
\[
\begin{array}{crl}
S&=&L_{G_{1}}+n_{2}I_{n_{1}}-\left(
  \begin{array}{cccccccccccccccc}
   0_{n_{1}\times n_{1}m_{2}}&-\mathbf1^{T}_{n_{2}}\otimes I_{n_{1}}\\
 \end{array}
  \right)
 \begin{array}{crl}
\left(
  \begin{array}{cccccccccccccccc}
    2I_{m_{2}}\otimes I_{n_{1}}& -R_{2}^{T}\otimes I_{n_{1}}\\
   -R_{2}\otimes I_{n_{1}} &A\\
 \end{array}
  \right)^{-1}
 \left(
  \begin{array}{cccccccccccccccc}
  0_{n_{1}m_{2}\times n_{1}}\\
  -\mathbf1_{n_{2}}\otimes I_{n_{1}}\\
 \end{array}
 \right),
\end{array}\\
\end{array}
\]
where $A=(D_{G_{2}}+I_{n_{2}})\otimes I_{n_{1}}$, then
\begin{eqnarray*}
S&=&L_{G_{1}}+n_{2}I_{n_{1}}-(\mathbf1^{T}_{n_{2}}((D_{G_{2}}+I_{n_{2}})I_{n_{2}}-\frac{1}{2}R_{2}R^{T}_{2})^{-1}\mathbf1_{n_{2}})\otimes I_{n_{1}})\\
&=&L_{G_{1}}+n_{2}I_{n_{1}}-(\mathbf1^{T}_{n_{2}}(\frac{1}{2}L_{G_{2}}+I_{n_{2}})^{-1}\mathbf1_{n_{2}})\otimes I_{n_{1}}\\
&=&L_{G_{1}}.
\end{eqnarray*}
By Lemma 2.4, we have $S^{\#}=L_{G_{1}}^{\#}$.
\vskip 0.1in
Next, according to Lemma 2.4, we calculate $-L^{-1}_{1}L_{2}S^{\#}$
and $-L^{-1}_{1}L_{2}S^{\#}L^{T}_{2}L^{-1}_{1}$.
\[
\begin{array}{crl}
-L^{-1}_{1}L_{2}S^{\#}
&=&-\left(
  \begin{array}{cccccccccccccccc}
    T_{1}& (R^{T}_{2}(L_{G_{2}}+2I_{n_{2}})^{-1})\otimes I_{n_{1}}\\
  ((L_{G_{2}}+2I_{n_{2}})^{-1}R_{2})\otimes I_{n_{1}} & 2(L_{G_{2}}+2I_{n_{2}})^{-1}\otimes I_{n_{1}}\\
 \end{array}
  \right)
 \begin{array}{crl}
 \left(
  \begin{array}{cccccccccccccccc}
    0_{n_{1}m_{2}\times n_{1}}\\
  -\mathbf1_{n_{2}}\otimes I_{n_{1}}\\
 \end{array}
 \right)S^{\#}
\end{array}\\
&=&\left(
  \begin{array}{cccccccccccccccc}
  (R^{T}_{2}(L_{G_{2}}+2I_{n_{2}})^{-1}\mathbf1_{n_{2}})\otimes I_{n_{1}}\\
  (2(L_{G_{2}}+2I_{n_{2}})^{-1}\mathbf1_{n_{2}})\otimes I_{n_{1}}\\
 \end{array}
 \right)S^{\#}.
\end{array}\\
\]
Note that $(L_{G_{2}}+2I_{n_{2}})\mathbf1_{n_{2}}=2\cdot\mathbf1_{n_{2}}$, then
$$R^{T}_{2}(L_{G_{2}}+2I_{n_{2}})^{-1}\mathbf1_{n_{2}}=\frac{1}{2}R^{T}_{2}\mathbf1_{n_{2}}=\mathbf1_{m_{2}},
~~~~~~~2(L_{G_{2}}+2I_{n_{2}})^{-1}\mathbf1_{n_{2}}=\mathbf1_{n_{2}},$$
For convenience, let $H=\mathbf1_{m_{2}}\otimes I_{n_{1}}, K=\mathbf1_{n_{2}}\otimes I_{n_{1}},$
then
\[
\begin{array}{crl}
-L^{-1}_{1}L_{2}S^{\#}
 &=&\left(
  \begin{array}{cccccccccccccccc}
  \mathbf1_{m_{2}}\otimes I_{n_{1}}\\
  \mathbf1_{n_{2}}\otimes I_{n_{1}}\\
 \end{array}
 \right)S^{\#}\\
 &=&\left(
  \begin{array}{cccccccccccccccc}
  HS^{\#}\\
  KS^{\#}\\
 \end{array}
 \right),\\
\end{array}
\]
and let $Q=(L_{G_{2}}+2I_{n_{2}})$, then
\[
\begin{array}{crl}
-L^{-1}_{1}L_{2}S^{\#}L^{T}_{2}L^{-1}_{1}
&=&\begin{array}{crl}
 \left(
  \begin{array}{cccccccccccccccc}
    -HS^{\#}\\
  -KS^{\#}\\
 \end{array}
 \right) \left(
  \begin{array}{cccccccccccccccc}
   0_{n_{1}\times n_{1}m_{2}}&-K^{T}\\
 \end{array}
 \right) \left(
  \begin{array}{cccccccccccccccc}
    T_{1}& (R^{T}_{2}Q^{-1})\otimes I_{n_{1}}\\
  (Q^{-1}R_{2})\otimes I_{n_{1}} & 2Q^{-1}\otimes I_{n_{1}}\\
 \end{array}
 \right)
\end{array}\\
&=&\left(
  \begin{array}{cccccccccccccccc}
-HS^{\#}\\
  -KS^{\#}\\
 \end{array}
 \right)\left(
  \begin{array}{cccccccccccccccc}
 (-\mathbf1_{n_{2}}^{T}Q^{-1}R_{2})\otimes I_{n_{1}} & -(2\cdot\mathbf1_{n_{2}}^{T}Q^{-1}\otimes I_{n_{1}}\\
 \end{array}
 \right).
\end{array}\\
\]
Note that $\mathbf1^{T}_{n_{2}}(L_{G_{2}}+2I_{n_{1}})=2\cdot\mathbf1^{T}_{n_{2}}$, then
$\mathbf1^{T}_{n_{2}}Q^{-1}R_{2}=\mathbf1^{T}_{n_{2}}(L_{G_{2}}+2I_{n_{1}})^{-1}R_{2}=\frac{1}{2}\mathbf1^{T}_{n_{2}}R_{2}=\mathbf1^{T}_{m_{2}}$,
so
\[
\begin{array}{crl}
-L^{-1}_{1}L_{2}S^{\#}L^{T}_{2}L^{-1}_{1}
&=&\begin{array}{crl}
 \left(
  \begin{array}{cccccccccccccccc}
    -HS^{\#}\\
  -KS^{\#}\\
 \end{array}
 \right) \left(
  \begin{array}{cccccccccccccccc}
   -\mathbf1^{T}_{m_{2}}\otimes I_{n_{1}}&-\mathbf1^{T}_{n_{2}}\otimes I_{n_{1}}\\
 \end{array}
 \right)
\end{array}\\
&=&\left(
  \begin{array}{cccccccccccccccc}
 HS^{\#}H^{T}& HS^{\#}K^{T}\\
 KS^{\#}H^{T}&KS^{\#}K^{T}\\
 \end{array}
 \right).
\end{array}\\
\]
Based on Lemma 2.3 and 2.4, the following matrix
\[
\begin{array}{crl}
N=\left(
  \begin{array}{cccccccccccccccc}
    T_{1}+HS^{\#}H^{T}&R^{T}_{2}Q^{-1}\otimes I_{n_{1}}+HS^{\#}K^{T}&HS^{\#}\\
Q^{-1}R_{2}\otimes I_{n_{1}}+KS^{\#}H^{T}& 2Q^{-1}\otimes I_{n_{1}}+KS^{\#}K^{T}&KS^{\#}\\
S^{\#}H&S^{\#}K&S^{\#}\\
 \end{array}
  \right)
\end{array}
~~~~~~~~~~~~~(3.1)
\]
is a symmetric $\{1\}$- inverse of $L(G_{1}\diamondsuit G_{2})$,
where $T_{1}=\frac{1}{2}(I_{m_{2}}+R^{T}_{2}(2I_{n_{2}}+L_{G_{2}})^{-1}R_{2})\otimes I_{n_{1}}$,
$Q=(L_{G_{2}}+2I_{n_{2}})$, $H=\mathbf1_{m_{2}}\otimes I_{n_{1}}, K=\mathbf1_{n_{2}}\otimes I_{n_{1}}$.

\vskip 0.1in
For any $i,j\in V(G_{1})$, by Lemma 2.1 and the Equation $(3.1)$, we have
\begin{eqnarray*}
r_{ij}(G_{1}\diamondsuit G_{2})&=&(L_{G_{1}}^{\#})_{ii}+(L_{G_{1}}^{\#})_{jj}-2(L_{G_{1}}^{\#})_{ij}.
\end{eqnarray*}

For any $i,j\in V(G_{2})$, by Lemma 2.1 and the Equation $(3.1)$, we have
\vskip 0.1in
~~~~~~~~~$r_{ij}(G_{1} \diamondsuit G_{2})=2((L_{G_{2}}+2I_{n_{2}})^{-1}\otimes I_{n_{1}})_{ii}+2((L_{G_{2}}+2I_{n_{2}})^{-1}\otimes I_{n_{1}})_{jj}$
\vskip 0.1in
$-2((L_{G_{2}}+2I_{n_{2}})^{-1}\otimes I_{n_{1}})_{ij}.$
\vskip 0.1in

For any $i\in V(G_{1})$, $j\in V(G_{2})$, by Lemma 2.1 and the Equation $(3.1)$, we have
\begin{eqnarray*}
r_{ij}(G_{1}\diamondsuit G_{2})&=&(L_{G_{1}}^{\#})_{ii}+2((L_{G_{2}}+2I_{n_{2}})^{-1}\otimes I_{n_{1}})_{jj}-(L_{G_{1}}^{\#})_{ij}.
\end{eqnarray*}

For any $i\in I(G_{2})$, $j\in V(G_{1})\cup V(G_{2})$, Let $u_{i}v_{i}\in E(G_{2})$ denote the edge corresponding
to $i$, By Lemma 2.5, we have
\begin{eqnarray*}
r_{ij}(G_{1}\diamondsuit G_{2})&=&\frac{1}{2}+\frac{1}{2}r_{u_{i}j}(G_{1}\diamondsuit G_{2})+\frac{1}{2}r_{v_{i}j}(G_{1}\diamondsuit G_{2})-\frac{1}{4}r_{u_{i}v_{i}}(G_{1}\diamondsuit G_{2}).
\end{eqnarray*}

For any $i,j\in I(G_{2})$, let $u_{i}v_{i}, u_{j}v_{j}\in E(G_{2})$ denote the edges corresponding
to $i,j$.  By Lemma 2.5, we have
\begin{eqnarray*}
r_{ij}(G_{1}\diamondsuit G_{2})
&=&\frac{1}{2}+\frac{1}{2}r_{u_{i}j}(G_{1}\diamondsuit G_{2})+\frac{1}{2}r_{v_{i}j}(G_{1}\diamondsuit G_{2})-\frac{1}{4}r_{u_{i}v_{i}}(G_{1}\diamondsuit G_{2})\\
&=& 1+\frac{1}{4}[r_{u_{i}u_{j}}(G_{1}\diamondsuit G_{2})+r_{u_{i}v_{j}}(G_{1}\diamondsuit G_{2})+r_{v_{i}u_{j}}(G_{1}\diamondsuit G_{2})\\
&&+r_{v_{i}v_{j}}(G_{1}\diamondsuit G_{2})-r_{u_{i}v_{i}}(G_{1}\diamondsuit G_{2})-r_{u_{j}v_{j}}(G_{1}\diamondsuit G_{2})].
\end{eqnarray*}
\hfill$\square$

When $G_{2}$ is a regular graph, we will give the formulae for resistance distance in $G_{1}\star G_{2}$ as follows.

{\bf Theorem 3.2} \  Let $G_{1}$ be an arbitrary graph on $n_{1}$ vertices, and
$G_{2}$ an $r_{2}$- regular graph on $n_{2}$ vertices and $m_{2}$ edges, then
the following hold:

(i) For any $i,j\in V(G_{1})$, we have
\begin{eqnarray*}
r_{ij}(G_{1}\star G_{2})&=&(L_{G_{1}}^{\#})_{ii}+(L_{G_{1}}^{\#})_{jj}-2(L_{G_{1}}^{\#})_{ij}.
\end{eqnarray*}

(ii) For any $i,j\in V(G_{2})$, we have
\begin{eqnarray*}
r_{ij}(G_{1}\star G_{2})&=&(3(r_{2}I_{n_{2}}+L_{G_{2}})\otimes I_{n_{1}})^{-1}_{ii}+(3(
r_{2}I_{n_{2}}+L_{G_{2}})\otimes I_{n_{1}})^{-1}_{jj}\\
&&-2(3(r_{2}I_{n_{2}}+L_{G_{2}})\otimes I_{n_{2}})^{-1}_{ij}.
\end{eqnarray*}

(iii) For any $i\in V(G_{1})$, $j\in V(G_{2})$, we have
\begin{eqnarray*}
r_{ij}(G_{1}\star G_{2})&=&(L_{G_{1}}^{\#})_{ii}+(3(r_{2}I_{n_{2}}+L(G_{2})\otimes I_{n_{1}})^{-1}_{jj}\\
&&-(L_{G_{1}}^{\#})_{ij}.
\end{eqnarray*}

(iv) For any $i\in I(G_{2})$, $j\in V(G_{1})\cup V(G_{2})$, let $u_{i}v_{i}\in E(G_{2})$ denote the edge corresponding
to $i$, $k\in V(G_{1})$, we have
\begin{eqnarray*}
r_{ij}(G_{1}\star G_{2})&=&\frac{1}{3}(1+r_{u_{i}j}(G_{1}\star G_{2})+r_{v_{i}j}(G_{1}\star G_{2})+r_{kj}(G_{1}\star G_{2}))\\
&&-\frac{1}{9}(r_{u_{i}v_{i}}(G_{1}\star G_{2})+r_{u_{i}k}(G_{1}\star G_{2})+r_{v_{i}k}(G_{1}\star G_{2})).
\end{eqnarray*}

(v) For any $i,j\in I(G_{2})$, let $u_{i}v_{i}, u_{j}v_{j}\in E(G_{2})$ denote the edges corresponding
to $i,j$, $k\in V(G_{1})$, we have
\begin{eqnarray*}
r_{ij}(G_{1}\star G_{2})&=&\frac{1}{3}(1+r_{ju_{i}}(G_{1}\star G_{2})+r_{jv_{i}}(G_{1}\star G_{2})+r_{jk}(G_{1}\star G_{2})\\
&&-\frac{1}{9}(r_{u_{i}v_{i}}(G_{1}\star G_{2})+r_{u_{i}k}(G_{1}\star G_{2})+r_{v_{i}k}(G_{1}\star G_{2}))\\
&=&\frac{1}{3}(1+r_{iu_{j}}(G_{1}\star G_{2})+r_{iv_{j}}(G_{1}\star G_{2})+r_{ik}(G_{1}\star G_{2})\\
&&-\frac{1}{9}(r_{u_{i}v_{i}}(G_{1}\star G_{2})+r_{u_{j}k}(G_{1}\star G_{2})+r_{v_{j}k}(G_{1}\star G_{2})).
\end{eqnarray*}

\vskip 0.1in
{\bf Proof.} \ Let $R_{2}$ be the incidence matrix of $G_{2}$. Then with
a proper labeling of vertices (see Ref.[13]),
the Laplacian matrix of $G_{1}\star G_{2}$ can be written as
\[
\begin{array}{crl}
 L(G_{1}\star G_{2})=\left(
  \begin{array}{cccccccccccccccc}
    3I_{m_{2}}\otimes I_{n_{1}}& -R_{2}^{T}\otimes I_{n_{1}}&-\mathbf1_{m_{2}}\otimes I_{n_{1}}\\
   -R_{2}\otimes I_{n_{1}} &r_{2}I_{n_{2}}\otimes I_{n_{1}}&0_{n_{1}n_{2}\times n_{1}}\\
   -\mathbf1^{T}_{m_{2}}\otimes I_{n_{1}}&0_{n_{1}\times n_{1}n_{2}}&L_{G_{1}}+m_{2}I_{n_{1}}\\
 \end{array}
  \right)
\end{array},
\]
where $0_{s,t}$ denotes the $s\times t$ matrix with all entries equal to zero. 
Let
\[
\begin{array}{crl}
 L_{1}=\left(
  \begin{array}{cccccccccccccccc}
   3I_{m_{2}}\otimes I_{n_{1}}& -R_{2}^{T}\otimes I_{n_{1}}\\
   -R_{2}\otimes I_{n_{1}} &r_{2}I_{n_{2}}\otimes I_{n_{1}}\\
 \end{array}
  \right)
\end{array},
\]
 $L_{2}=\left(
  \begin{array}{cc}
    -\mathbf1_{m_{2}}\otimes I_{n_{1}}\\
   0_{n_{1}n_{2}\times n_{1}} \\
  \end{array}
\right)$, $L_{2}^{T}=\left(
  \begin{array}{cc}
     -\mathbf1^{T}_{m_{2}}\otimes I_{n_{1}}&0_{n_{1}\times n_{1}n_{2}}\\
  \end{array}
\right)$, and $L_{3}=L_{G_{1}}+m_{2}I_{n_{1}}$.

\vskip 0.1in
We begin with the calculation about $L_{1}^{-1}$. By Lemma 2.3, we have
\begin{eqnarray*}
R_{2}&=&[r_{2}I_{n_{2}}\otimes I_{n_{1}}-(R_{2}\otimes I_{n_{1}})(3I_{m_{2}}\otimes I_{n_{1}})^{-1}(R_{2}^{T}\otimes I_{n_{1}})]^{-1}\\
&=&[r_{2}I_{n_{2}}\otimes I_{n_{1}}-\frac{1}{3}R_{2}R^{T}_{2}\otimes I_{n_{1}}]^{-1}\\
&=&[r_{2}I_{n_{2}}-\frac{1}{3}(r_{2}I_{n_{2}}+A_{G_{2}})]^{-1}\otimes I_{n_{1}}\\
&=&3(L_{G_{2}}+r_{2}I_{n_{2}})^{-1}\otimes I_{n_{1}}.
\end{eqnarray*}
One can obtain that
\[
\begin{array}{crl}
L_{1}^{-1}=\left(
  \begin{array}{cccccccccccccccc}
    T_{2}& (R^{T}_{2}(L_{G_{2}}+r_{2}I_{n_{2}})^{-1})\otimes I_{n_{1}}\\
  ((L_{G_{2}}+r_{2}I_{n_{2}})^{-1}R_{2})\otimes I_{n_{1}} & 3(L_{G_{2}}+r_{2}I_{n_{2}})^{-1}\otimes I_{n_{1}}\\
 \end{array}
  \right),
\end{array}
\]
where $T_{2}=\frac{1}{3}(I_{m_{2}}+R^{T}_{2}(r_{2}I_{n_{2}}+L_{G_{2}})^{-1}R_{2})\otimes I_{n_{1}}$.
\vskip 0.1in
Now we are ready to calculate $S$.
\vskip 0.1in

Let
\[
\begin{array}{crl}
S=L_{G_{1}}+m_{2}I_{n_{1}}-\left(
  \begin{array}{cccccccccccccccc}
   -\mathbf1^{T}_{m_{2}}\otimes I_{n_{1}}&0_{n_{1}\times n_{1}n_{2}}\\
 \end{array}
  \right)
 \begin{array}{crl}
\left(
  \begin{array}{cccccccccccccccc}
    3I_{m_{2}}\otimes I_{n_{1}}& -R_{2}^{T}\otimes I_{n_{1}}\\
   -R_{2}\otimes I_{n_{1}} &r_{2}I_{n_{2}}\otimes I_{n_{1}}&\\
 \end{array}
  \right)^{-1}
  \left(
  \begin{array}{cccccccccccccccc}
  -\mathbf1_{m_{2}}\otimes I_{n_{1}}\\
   0_{n_{1}n_{2}\times n_{1}}\\
 \end{array}
 \right),
\end{array}
\end{array}
\]
then
\begin{eqnarray*}
S&=&L_{G_{1}}+m_{2}I_{n_{1}}-(-\mathbf1^{T}_{m_{2}}\otimes I_{n_{1}})((3I_{m_{2}}-\frac{1}{r_{2}}R^{T}_{2}R_{2})^{-1}\otimes I_{n_{1}})(-\mathbf1_{m_{2}}\otimes I_{n_{1}})\\
&=&L_{G_{1}}+m_{2}I_{n_{1}}-(r_{2}\mathbf1^{T}_{m_{2}}(L_{l(G_{2})}+r_{2}I_{m_{2}})^{-1}\mathbf1_{m_{2}})\otimes I_{n_{1}}\\
&=&L_{G_{1}}.
\end{eqnarray*}
By Lemma 2.4, we have $S^{\#}=L_{G_{1}}^{\#}$.
\vskip 0.1in
Similiarly, according to Lemma 2.4, we calculate $-L^{-1}_{1}L_{2}S^{\#}$
and $-L^{-1}_{1}L_{2}S^{\#}L^{T}_{2}L^{-1}_{1}$.
\[
\begin{array}{crl}
-L^{-1}_{1}L_{2}S^{\#}
&=&\left(
  \begin{array}{cccccccccccccccc}
    T_{2}& (R^{T}_{2}(L_{G_{2}}+r_{2}I_{n_{2}})^{-1})\otimes I_{n_{1}}\\
  ((L_{G_{2}}+r_{2}I_{n_{2}})^{-1}R_{2})\otimes I_{n_{1}} & 3(L_{G_{2}}+r_{2}I_{n_{2}})^{-1}\otimes I_{n_{1}}\\
 \end{array}
  \right)
 \begin{array}{crl}
 \left(
  \begin{array}{cccccccccccccccc}
   -\mathbf1_{m_{2}}\otimes I_{n_{1}}\\
   0_{n_{1}n_{2}\times n_{1}}\\
 \end{array}
 \right)S^{\#}
\end{array}\\
&=&\left(
  \begin{array}{cccccccccccccccc}
  -T_{2}(\mathbf1_{m_{2}}\otimes I_{n_{1}})\\
  -((L_{G_{2}}+r_{2}I_{n_{2}})^{-1}R_{2})\otimes I_{n_{1}})(\mathbf1_{m_{2}}\otimes I_{n_{1}})\\
 \end{array}
 \right)S^{\#}.
\end{array}\\
\]

Note that $(L_{G_{2}}+r_{2}I_{n_{2}})\mathbf1_{n_{2}}=r_{2}\mathbf1_{n_{2}}$, then
\begin{eqnarray*}
T_{2}(1_{m_{2}}\otimes I_{n_{1}})&=&\frac{1}{3}((I_{m_{2}}+R^{T}_{2}(r_{2}I_{n_{2}}+L_{G_{2}})^{-1}R_{2})\otimes I_{n_{1}})(\mathbf1_{m_{2}}\otimes I_{n_{1}})\\
&=&\frac{1}{3}(\mathbf1_{m_{2}}+R_{2}^{T}(r_{2}I_{n_{2}}+L_{G_{2}})^{-1}R_{2}\mathbf1_{m_{2}})\otimes I_{n_{1}}\\
&=&\frac{1}{3}(\mathbf1_{m_{2}}+R_{2}^{T}(r_{2}I_{n_{2}}+L_{G_{2}})^{-1}r_{2}\mathbf1_{n_{2}})\otimes I_{n_{1}}\\
&=&\frac{1}{3}(\mathbf1_{m_{2}}+R^{T}_{2}\mathbf1_{n_{2}})\otimes I_{n_{1}}=\mathbf1_{m_{2}}\otimes I_{n_{1}}.
\end{eqnarray*}
and $((r_{2}I_{n_{2}}+L_{G_{2}})^{-1}R_{2})\otimes I_{n_{1}})(\mathbf1_{m_{2}}\otimes I_{n_{1}})=((r_{2}I_{n_{2}}+L_{G_{2}})^{-1}r_{2}\mathbf1_{n_{2}})
\otimes I_{n_{1}}=\mathbf1_{n_{2}}\otimes I_{n_{1}}$, so
\[
\begin{array}{crl}
-L^{-1}_{1}L_{2}S^{\#}
 &=&\left(
  \begin{array}{cccccccccccccccc}
  -(\mathbf1_{m_{2}}\otimes I_{n_{1}})\\
  -(\mathbf1_{n_{2}}\otimes I_{n_{1}})\\
 \end{array}
 \right)S^{\#}
 =\left(
  \begin{array}{cccccccccccccccc}
-(\mathbf1_{m_{2}}\otimes I_{n_{1}})S^{\#}\\
  -(\mathbf1_{n_{2}}\otimes I_{n_{1}})S^{\#}\\
 \end{array}
 \right).\\
\end{array}
\]
Let $C=L_{G_{2}}+r_{2}I_{n_{2}}$, $H=\mathbf1_{m_{2}}\otimes I_{n_{1}}, K=\mathbf1_{n_{2}}\otimes I_{n_{1}}$,
\[
\begin{array}{crl}
-L^{-1}_{1}L_{2}S^{\#}L^{T}_{2}L^{-1}_{1}
&=&\begin{array}{crl}
 \left(
  \begin{array}{cccccccccccccccc}
   -HS^{\#}\\
 -KS^{\#}\\
 \end{array}
 \right) \left(
  \begin{array}{cccccccccccccccc}
   -H^{T}&0_{n_{1}\times n_{1}n_{2}}\\
 \end{array}
 \right) \left(
  \begin{array}{cccccccccccccccc}
    T_{2}& R^{T}_{2}C^{-1}\otimes I_{n_{1}}\\
  C^{-1}R_{2}\otimes I_{n_{1}} & 3C^{-1}\otimes I_{n_{1}}\\
 \end{array}
 \right)
\end{array}\\
&=&\left(
  \begin{array}{cccccccccccccccc}
 -HS^{\#}\\
  -KS^{\#}\\
 \end{array}
 \right)\left(
  \begin{array}{cccccccccccccccc}
-H^{T}T_{2}&-H^{T}(R^{T}_{2}C^{-1}\otimes I_{n_{1}})\\
 \end{array}
 \right).
\end{array}\\
\]
Note that $H^{T}T_{2}=(\mathbf1_{m_{2}}\otimes I_{n_{1}})^{T}T_{2}=(T_{2}(\mathbf1_{m_{2}}\otimes I_{n_{1}}))^{T}=H^{T}$, then
\[
\begin{array}{crl}
-L^{-1}_{1}L_{2}S^{\#}L^{T}_{2}L^{-1}_{1}
&=&\begin{array}{crl}
 \left(
  \begin{array}{cccccccccccccccc}
   -HS^{\#}\\
  -KS^{\#}\\
 \end{array}
 \right) \left(
  \begin{array}{cccccccccccccccc}
 -H^{T}&K^{T}\\
 \end{array}
 \right)
\end{array}
=\left(
  \begin{array}{cccccccccccccccc}
 HS^{\#}H^{T}&HS^{\#}K^{T}\\
KS^{\#}H^{T}&KS^{\#}K^{T}\\
 \end{array}
 \right).
\end{array}\\
\]
Based on Lemma 2.3 and 2.4, the following matrix
\[
\begin{array}{crl}
N=\left(
  \begin{array}{cccccccccccccccc}
    T_{2}+HS^{\#}H^{T}&R^{T}_{2}C^{-1}\otimes I_{n_{1}}+HS^{\#}K^{T}&HS^{\#}\\
C^{-1}R_{2}\otimes I_{n_{1}}+KS^{\#}H^{T}& 3C^{-1}\otimes I_{n_{1}}+KS^{\#}K^{T}&KS^{\#}\\
S^{\#}H&S^{\#}K&S^{\#}\\
 \end{array}
  \right)
\end{array}
~~~~~~~~~~~~~(3.2)
\]
is a symmetric $\{1\}$- inverse of $L(G_{1}\star G_{2})$,
where $T_{2}=\frac{1}{3}(I_{m_{2}}+R^{T}_{2}(r_{2}I_{n_{2}}+L_{G_{2}})^{-1}R_{2})\otimes I_{n_{1}}$,
$C=(L_{G_{2}}+r_{2}I_{n_{2}})$, $H=\mathbf1_{m_{2}}\otimes I_{n_{1}}, K=\mathbf1_{n_{2}}\otimes I_{n_{1}}$.
\vskip 0.1in
For any $i,j\in V(G_{1})$, by Lemma 2.1 and the Equation $(3.2)$, we have
\begin{eqnarray*}
r_{ij}(G_{1}\star G_{2})&=&(L_{G_{1}}^{\#})_{ii}+(L_{G_{1}}^{\#})_{jj}-2(L_{G_{1}}^{\#})_{ij}.
\end{eqnarray*}

For any $i,j\in V(G_{2})$, by Lemma 2.1 and the Equation $(3.2)$, we have
\begin{eqnarray*}
r_{ij}(G_{1}\star G_{2})&=&(3(r_{2}I_{n_{2}}+L_{G_{2}})\otimes I_{n_{1}})^{-1}_{ii}+(3(
r_{2}I_{n_{2}}+L_{G_{2}})\otimes I_{n_{1}})^{-1}_{jj}\\
&&-2(3(r_{2}I_{n_{2}}+L_{G_{2}})\otimes I_{n_{1}})^{-1}_{ij}.
\end{eqnarray*}

For any $i\in V(G_{1})$, $j\in V(G_{2})$, by Lemma 2.1 and the Equation $(3.2)$, we have
\begin{eqnarray*}
r_{ij}(G_{1}\star G_{2})&=&(L_{G_{1}}^{\#})_{ii}+(3(r_{2}I_{n_{2}}+L(G_{2})\otimes I_{n_{1}})^{-1}_{jj}\\
&&-(L_{G_{1}}^{\#})_{ij}.
\end{eqnarray*}

For any $i\in I(G_{2})$, $j\in V(G_{1})\cup V(G_{2})$, let $u_{i}v_{i}\in E(G_{2})$ denote the edge corresponding
to $i$, $k\in V(G_{1})$, by Lemma 2.5, we have
\begin{eqnarray*}
r_{ij}(G_{1}\star G_{2})&=&\frac{1}{3}(1+r_{u_{i}j}(G_{1}\star G_{2})+r_{v_{i}j}(G_{1}\star G_{2})+r_{kj}(G_{1}\star G_{2}))\\
&&-\frac{1}{9}(r_{u_{i}v_{i}}(G_{1}\star G_{2})+r_{u_{i}k}(G_{1}\star G_{2})+r_{v_{i}k}(G_{1}\star G_{2})).
\end{eqnarray*}

For any $i,j\in I(G_{2})$, let $u_{i}v_{i}, u_{j}v_{j}\in E(G_{2})$ denote the edges corresponding
to $i,j$, $k\in V(G_{1})$, by Lemma 2.5, we have
\begin{eqnarray*}
r_{ij}(G_{1}\star G_{2})&=&\frac{1}{3}(1+r_{ju_{i}}(G_{1}\star G_{2})+r_{jv_{i}}(G_{1}\star G_{2})+r_{jk}(G_{1}\star G_{2})\\
&&-\frac{1}{9}(r_{u_{i}v_{i}}(G_{1}\star G_{2})+r_{u_{i}k}(G_{1}\star G_{2})+r_{v_{i}k}(G_{1}\star G_{2}))\\
&=&\frac{1}{3}(1+r_{iu_{j}}(G_{1}\star G_{2})+r_{iv_{j}}(G_{1}\star G_{2})+r_{ik}(G_{1}\star G_{2})\\
&&-\frac{1}{9}(r_{u_{i}v_{i}}(G_{1}\star G_{2})+r_{u_{j}k}(G_{1}\star G_{2})+r_{v_{j}k}(G_{1}\star G_{2})).
\end{eqnarray*}
\hfill$\square$

\section{Kirchhoff index in corona-vertex and
corona-edge of subdivision graph}

{\bf Theorem 4.1} \  Let $G_{1}$ be a graph on $n_{1}$ vertices and $m_{1}$
edges and $G_{2}$ be a graph on $n_{2}$ vertices and $m_{2}$ edges.
Then
\begin{eqnarray*}
Kf(G_{1}\diamondsuit G_{2})
&=&n_{1}(1+n_{2}+m_{2})[\frac{n_{1}m_{2}}{2}+\frac{n_{1}}{2}(tr(Q^{-1}A_{G_{2}})+tr(Q^{-1}D_{G_{2}}))+2n_{1}\sum_{i=1}^{n_{2}}\frac{1}{\mu_{i}(G_{2})+2}\\
&&+\frac{m_{2}+n_{2}+1}{n_{1}}Kf(G_{1})]-\frac{n_{1}}{2}(\pi^{T}Q^{-1}\pi)-\frac{5n_{1}m_{2}+2n_{1}n_{2}}{2},
\end{eqnarray*}
where $Q=L_{G_{2}}+2I_{n_{2}}$, $D_{G_{2}}=diag(d_{1}, d_{2},\cdots, d_{n_{2}})$ and
$\pi=(d_{1}, d_{2},\cdots, d_{n_{2}})^{T}$.
\vskip 0.1in
{\bf Proof.} \ Let $L^{(1)}_{G_{1}\diamondsuit G_{2}}$ be the symmetric $\{1\}$-inverse of $L_{G_{1}\diamondsuit G_{2}}$.
Then
\begin{eqnarray*}
tr(L^{(1)}_{G_{1}\diamondsuit G_{2}})
&=&tr(T_{1}+HS^{\#}H^{T})+tr(2Q^{-1}\otimes I_{n_{1}}+KS^{\#}K^{T})+tr(S^{\#}),
\end{eqnarray*}
where $T_{1}=\frac{1}{2}(I_{m_{2}}+R^{T}_{2}(2I_{n_{2}}+L_{G_{2}})^{-1}R_{2})\otimes I_{n_{1}}$, $Q=(L_{G_{2}}+2I_{n_{2}})$, $S^{\#}=L^{\#}_{G_{1}}$, $H=\mathbf1_{m_{2}}\otimes I_{n_{1}}, K=\mathbf1_{n_{2}}\otimes I_{n_{1}}$,
and $R_{2}$ is the incidence matrix of $G_{2}$.
\vskip 0.1in
By Lemma 2.7, $HS^{\#}H^{T}=j_{m_{2}\times m_{2}}\otimes S^{\#}$, $KS^{\#}K^{T}=j_{n_{2}\times n_{2}}\otimes S^{\#}$,
then $$tr(HS^{\#}H^{T})=m_{2}tr(S^{\#}),~~~~tr(KS^{\#}K^{T})=n_{2}tr(S^{\#}).$$
Note that the eigenvalues of $Q$ are $\mu_{1}(G_{2})+2, ...,\mu_{n}(G_{2})+2$, then
$$tr(Q^{-1})=\sum_{i=1}^{n_{2}}\frac{1}{\mu_{i}(G_{2})+2}.$$
Recall that $R_{2}R_{2}^{T}=D_{G_{2}}+A_{G_{2}}$, then
$$tr(R_{2}^{T}Q^{-1}R_{2})=tr(Q^{-1}R_{2}R_{2}^{T})=tr(Q^{-1}A_{G_{2}})+tr(Q^{-1}D_{G_{2}}).$$

Thus, substituting $tr(R_{2}^{T}Q^{-1}R_{2})$, $tr(Q^{-1})$ into $tr(L^{(1)}_{G_{1}\diamondsuit G_{2}})$, we have
\begin{eqnarray*}
tr(L^{(1)}_{G_{1}\diamondsuit G_{2}})
&=&tr(T_{1})+tr(j_{m_{2}\times m_{2}}\otimes S^{\#})+2tr(Q^{-1}\otimes I_{n_{1}})+tr(j_{n_{2}\times n_{2}}\otimes S^{\#})+tr(S^{\#})\\
&=&\frac{n_{1}m_{2}}{2}+\frac{n_{1}}{2}tr(R_{2}^{T}Q^{-1}R_{2})+m_{2}tr(S^{\#})+2n_{1}tr(Q^{-1})+n_{2}tr(S^{\#})+tr(S^{\#})\\
&=&\frac{n_{1}m_{2}}{2}+\frac{n_{1}}{2}(tr(Q^{-1}A_{G_{2}})+tr(Q^{-1}D_{G_{2}}))+m_{2}tr(S^{\#})\\
&&+2n_{1}\sum_{i=1}^{n_{2}}\frac{1}{\mu_{i}(G_{2})+2}+n_{2}tr(S^{\#})+tr(S^{\#})\\
&=&\frac{n_{1}m_{2}}{2}+\frac{n_{1}}{2}(tr(Q^{-1}A_{G_{2}})+tr(Q^{-1}D_{G_{2}}))+2n_{1}\sum_{i=1}^{n_{2}}\frac{1}{\mu_{i}(G_{2})+2}\\
&&+\frac{m_{2}+n_{2}+1}{n_{1}}Kf(G_{1}).
\end{eqnarray*}
Next, we calculate the $\mathbf{1}^{T}(L^{(1)}_{G_{1}\diamondsuit G_{2}})\mathbf{1}$.
Since $L_{G}^{\#}\mathbf{1}=0$, by Lemma 2.7, then
\begin{eqnarray*}
\mathbf{1}^{T}(L^{(1)}_{G_{1}\diamondsuit G_{2}})\mathbf{1}
&=&\mathbf{1}^{T}T_{1}\mathbf{1}+\mathbf1^{T}(2Q^{-1}\otimes I_{n_{1}})\mathbf{1}+\mathbf{1}^{T}((R^{T}_{2}Q^{-1})\otimes I_{n_{1}})\mathbf{1}+\mathbf{1}^{T}((Q^{-1}R_{2}^{T})\otimes I_{n_{1}})\mathbf{1}.
\end{eqnarray*}
Note that $R_{2}\mathbf{1}=\pi$, where $\pi=(d_{1}, d_{2},\cdots, d_{n_{2}})^{T}$,
then $$\mathbf{1}^{T}((R_{2}^{T}Q^{-1}R_{2})\otimes I_{n_{1}})\mathbf{1}=n_{1}(\pi^{T}Q^{-1}\pi).$$
Since $Q\mathbf1_{n_{2}}=2\cdot\mathbf1_{n_{2}}$, then
$\mathbf{1}^{T}(2Q^{-1}\otimes I_{n_{1}})\mathbf{1}=2n_{1}(\mathbf{1}^{T}Q^{-1}\mathbf{1})=n_{1}n_{2}$ and
$$\mathbf{1}^{T}((Q^{-1}R_{2})\otimes I_{n_{1}})\mathbf{1}=n_{1}\mathbf{1}^{T}(Q^{-1}R_{2})\mathbf{1}
=\frac{1}{2}\mathbf{1}^{T}R_{2}\mathbf{1}
=\frac{1}{2}n_{1}\sum_{i=1}^{n_{2}}d_{i}=n_{1}m_{2},$$
Similarly, $\mathbf{1}^{T}((R_{2}^{T}Q^{-1})\otimes I_{n_{1}})\mathbf{1}=(\mathbf{1}^{T}((Q^{-1}R_{2})\otimes I_{n_{1}})\mathbf{1})^{T}=n_{1}m_{2}$.
\vskip 0.1in
Substituting $\mathbf{1}^{T}((R_{2}^{T}Q^{-1}R_{2})\otimes I_{n_{1}})\mathbf{1}$, $\mathbf{1}^{T}((Q^{-1}R_{2})\otimes I_{n_{1}})\mathbf{1}$, $\mathbf{1}^{T}((R_{2}^{T}Q^{-1})\otimes I_{n_{1}})\mathbf{1}$
into $\mathbf{1}^{T}(L^{(1)}_{G_{1}\diamondsuit G_{2}})\mathbf{1}$, we get
\begin{eqnarray*}
\mathbf{1}^{T}(L^{(1)}_{G_{1}\diamondsuit G_{2}})\mathbf{1}
&=&\frac{n_{1}m_{2}}{2}+\frac{n_{1}}{2}(\pi^{T}Q^{-1}\pi)
+2m_{2}n_{1}+n_{1}n_{2}.
\end{eqnarray*}
Lemma 2.6 implies that
$$Kf(L^{(1)}_{G_{1}\diamondsuit G_{2}})=n_{1}(1+n_{2}+m_{2})tr(L^{(1)}_{G_{1}\diamondsuit G_{2}})-\mathbf{1}^{T}(L^{(1)}_{G_{1}\diamondsuit G_{2}})\mathbf{1}.$$
Then plugging $tr(L^{(1)}_{G_{1}\diamondsuit G_{2}})$ and $\mathbf{1}^{T}(L^{(1)}_{G_{1}\diamondsuit G_{2}})\mathbf1$
into the equation above, we obtain the required result.\hfill$\square$

\vskip 0.1in
{\bf Corollary 4.2} \ Let $G_{1}$ be a graph on $n_{1}$ vertices and $m_{1}$ edges
and $G_{2}$ an $r_{2}$- regular graph on $n_{2}$ vertices and $m_{2}$ edges. Then
\begin{eqnarray*}
Kf(G_{1}\diamondsuit G_{2})
&=&n_{1}(1+n_{2}+m_{2})[\frac{n_{1}m_{2}}{2}+\frac{n_{1}}{2}(2r_{2}\sum_{i=1}^{n_{2}}\frac{1}{\mu_{i}(G_{2})+2}
-\sum_{i=1}^{n_{2}}\frac{\mu_{i}(G_{2})}{\mu_{i}(G_{2})+2})+\\
&&2n_{1}\sum_{i=1}^{n_{2}}\frac{1}{\mu_{i}(G_{2})+2}
+\frac{m_{2}+n_{2}+1}{n_{1}}Kf(G_{1})]-\frac{n_{1}n_{2}r^{2}_{2}}{4}-\frac{5n_{1}m_{2}+2n_{1}n_{2}}{2}.
\end{eqnarray*}
\vskip 0.1in
{\bf Proof.} \ Since $tr(Q^{-1}D_{G_{2}})=r_{2}tr(Q^{-1})=r_{2}\sum_{i=1}^{n_{2}}\frac{1}{\mu_{i}(G_{2})+2}$,
$tr(Q^{-1}A_{G_{2}})=r_{2}tr(Q^{-1})
-tr(Q^{-1}L_{G_{2}})=r_{2}\sum_{i=1}^{n_{2}}\frac{1}{\mu_{i}(G_{2})+2}
-\sum_{i=1}^{n_{2}}\frac{\mu_{i}(G_{2})}{\mu_{i}(G_{2})+2}$, $\pi^{T}Q^{-1}\pi=r^{2}\mathbf{1}^{T}Q^{-1}\mathbf{1}=
\frac{n_{2}r^{2}_{2}}{2}$, then the required result is obtained by plugging
$tr(Q^{-1}D_{G_{2}})$, $tr(Q^{-1}A_{G_{2}})$ and $\pi^{T}Q^{-1}\pi$ into
Theorem 4.1.
\hfill$\square$

{\bf Theorem 4.3} \ Let $G_{1}$ be an arbitrary graph on $n_{1}$ vertices, and
$G_{2}$ an $r_{2}$- regular graph on $n_{2}$ vertices and $m_{2}$ edges. Then
\begin{eqnarray*}
Kf(G_{1}\star G_{2})
&=&n_{1}(1+n_{2}+m_{2})[\frac{n_{1}m_{2}}{3}+\frac{n_{1}}{3}(tr(C^{-1}A_{G_{2}})+r_{2}\sum_{i=1}^{n_{2}}\frac{1}{\mu_{i}(G_{2})+r_{2}})+\\
&&3n_{1}\sum_{i=1}^{n_{2}}\frac{1}{\mu_{i}(G_{2})+r_{2}}+\frac{m_{2}+n_{2}+1}{n_{1}}Kf(G_{1})]-\frac{n_{1}m_{2}r_{2}+n_{1}n_{2}(r_{2}+3)^{2}}{3r_{2}},
\end{eqnarray*}
where $C=L_{G_{2}}+r_{2}I_{n_{2}}$.
\vskip 0.1in
{\bf Proof.} \ The proof is similar to that of Theorem 4.1 and hence we omit details.\hfill$\square$

\vskip 0.1in \noindent{\bf Acknowledgment:} This work was
supported by the National Natural Science Foundation of China
(Nos. 11461020, 11561042, 11601006 and 11471016), the Youth
Foundation of Hexi University in Gansu Province (No.QN2013-07) and
the Teaching Reform Project of Hexi University (No.
HXXXJY-2014-011).

\end{document}